\documentclass[10pt]{article}
\usepackage{txfonts}
\usepackage{graphicx}
\usepackage{fancyheadings}

\begin{document}

\pagenumbering{arabic}
\pagestyle{fancyplain}
\cfoot{\rm\thepage}
\lhead[18th World IMACS / MODSIM Congress, Cairns, Australia 13-17 July 2009
http://mssanz.org.au/modsim09]{18th World IMACS / MODSIM Congress, Cairns, Australia 13-17 July 2009\\
http://mssanz.org.au/modsim09}
\rhead[]{}

\begin{center}
{\Large{\bf{CHASSIS - Inverse Modelling of Relaxed Dynamical Systems\\}}}
{\large{Dalia Chakrabarty \\
School of Physics $\&$ Astronomy, University of Nottingham, 
Nottingham NG7 2RD, U.K.\\
email: dalia.chakrabarty@nottingham.ac.uk}}
\end{center}
\date{}


\noindent
The state of a non-relativistic gravitational dynamical system is known
at any time $t$ if the dynamical rule, i.e. Newton's equations of
motion, can be solved; this requires specification of the
gravitational potential. The evolution of a bunch of phase space
coordinates ${\bf w}$ is deterministic, though generally
non-linear. We discuss the novel Bayesian non-parametric algorithm
CHASSIS that gives phase space $pdf$ $f({\bf w})$ and potential
$\Phi({\bf x})$ of a relaxed gravitational system. CHASSIS is
undemanding in terms of input requirements in that it is viable given
incomplete, single-component velocity information of system
members. Here ${\bf x}$ is the 3-D spatial coordinate and ${\bf
w}={\bf x+v}$ where ${\bf v}$ is the 3-D velocity vector. CHASSIS
works with a 2-integral $f=f(E, L)$ where energy $E=\Phi + v^2/2, \:
v^2 = \sum_{i=1}^{3}{v_i^2}$ and the angular momentum is $L = |{\bf
r}\times{\bf v}|$, where ${\bf r}$ is the spherical spatial
vector. Also, we assume spherical symmetry. CHASSIS obtains the
$f(\cdot)$ from which the kinematic data is most likely to have been
drawn, in the best choice for $\Phi(\cdot)$, using an MCMC optimiser
(Metropolis-Hastings). The likelihood function ${\cal{L}}$ is defined
in terms of the projections of $f(\cdot)$ into the space of
observables and the maximum in ${\cal{L}}$ is sought by the optimiser.
The recovered solutions can be susceptible to large uncertainties
given the dimensionality of the domain of the unconstrained $f(\cdot)$
and the typically small, observed velocity samples in distant
astrophysical systems. This scenario is tackled by assuming $f=f(E)$,
i.e. we assume the phase space to be isotropic. However, this
simplifying assumption of isotropy is addressed by undertaking a
Bayesian test of significance that is developed to be used in the
non-parametric context.

A test based on the $p$-value estimates of the goodness of isotropy in
the data was previously undertaken (Chakrabarty $\&$ Saha
2001). However, $p$-values are sensitive to sample sizes and obfuscate
interpretation of analyses of differently sized kinematic
samples. Thus, a Bayesian formalism is a better alternative, eg. Fully
Bayesian Significance Test or FBST (Pereira $\&$ Stern 1999, Pereira,
Stern $\&$ Wechsler 2008). The null hypothesis that we aim to test, is
that the data are drawn from an isotropic $f(\cdot)$, i.e. $H_0:
{\hat{f}}=\Psi[E(v^2/2 + \Phi(r))]$ where the data are drawn from
$\hat{f}$ and $\Psi$ is some function: $\Psi > 0$ for $E < 0$ and
$\Psi=0$ otherwise. Within FBST, the evidence value ($ev$) in favour
of $H_0$ is obtained by first numerically spotting the most likely
configuration ($\theta^*$) that is compatible with $H_0$ and then
finding the volume of the tangential set $T$ by numerical
integration. Here $T$ is the set of all configurations with posterior
probability in excess of that of $\theta^*$. We have developed the
implementation of this scheme in the non-parametric context; in
CHASSIS, the configurations are $f(E)$ and $\Phi(r)$. To have
configurations obeying $H_0$, we perform sampling from the
$f(E)-\Phi(r)$ pair identified upon convergence of a run of
CHASSIS. From this sampling, the resulting $f(E)-\Phi(r)$
configuration corresponding to the highest ${\cal{L}}$ is compared to
all the other $f(E)-\Phi(r)$ pairs, in order to obtain a measure of
the volume of $T$.

We discuss two distinct applications of the isotropic version of
CHASSIS. In one, 2 distinct kinematic data sets of 2 distinct types of
members of an example galaxy are analysed by CHASSIS under the
assumption of isotropy. The $f(\cdot)$ and $\Phi(\cdot)$ recovered
from runs done with the two data sets are identified as
distinct. Given that the same galaxy cannot be described by two
different gravitational potentials, the risk involved in the very
method of extracting the galactic potential from kinematic data of
individual galactic members is demonstrated here for the first
time. The goodness of the assumption of isotropy, given the 2 data
sets, is quantified by our Bayesian test of hypothesis.

In the second application, it is shown that once the amount of
gravitational matter inside a fiduciary radius is pinned down from
independent measurements, the recovered $\Phi(\cdot)$ is 
unique, irrespective of the toy $f(\cdot)$, from which kinematic data
samples are drawn as input for CHASSIS, as long as velocity dispersion
values are measured at 1 or more different locations in the system. This
consistent nature of $\Phi(\cdot)$ is arrived at, notwithstanding
varied forms for the assumed toy $f(\cdot)$, including isotropic as
well as $E\:\&\:L$ dependent forms.\\

\clearpage

\lhead[Dalia Chakrabarty: CHASSIS]{Dalia Chakrabarty: CHASSIS}
\rhead[]{}

\section{Introduction}
\noindent
The complete characterisation of a gravitationally bound,
non-relativistic dynamical system can be undertaken with the help of
$f({\bf{w}}, t)$ - the $pdf$ of phase space $W$ - and the
gravitational potential $\Phi(\bf{x})$; here $t$ is time and ${\bf{w=
x+v}}$, where ${\bf{x}}$ represents the spatial coordinates and the
velocity vector is ${\bf{v}}={\bf{\dot{x}}}$. A sample of phase space
coordinates can be drawn from $f(\cdot)$ and allowed to evolve in
$\Phi(\cdot)$, in accordance to Newton's laws. In this way, the
evolution of the system is deterministic at any time $t$, though
non-linear in general. Hence, we aim to estimate $f(\cdot)$ and
$\Phi(\cdot)$; we focus on the characterisation of astrophysical
systems in this paper. A related aim is to derive the distribution of
the total gravitating matter from the estimated $\Phi(\cdot)$, while
keeping in mind that such $total$ mass is accounted for only partly by
luminous matter while the greater fraction is dark matter in these
systems.

Conventionally, 
\begin{itemize}
\item $f(\cdot)$ and $\Phi(\cdot)$ are almost always
parametrically described. However, given that astrophysical systems
such as galaxies are more likely to manifest complexity in their
dynamics than otherwise, any smooth parametric description of such
systems is erroneous. 
\item mass determination is typically pursued via observed photometric
  or luminous information though no functional dependence of the total
  (luminous+dark) matter content on such measurements 
  exist.
\item inhomogeneities in the measurement
errors notwithstanding, goodness-of-fit parameters are often invoked to
seek the solution. For galaxies, measurements are typically noisy and
such goodness-of-fit parameters can be artificially inflated (Bissantz $\&$ Munk, 2001). 
\end{itemize}
All these issues suggest a better - preferably, a nonparametric -
route to $f(\cdot)$ and $\Phi(\cdot)$ determination. It is such a {\it
novel, data driven characterisation of real 
(as
distinguished from simulated) 
astrophysical systems that CHASSIS
offers, using the few kinematic measurements that are typically
available}. 
Importantly, we test the chief assumption of our algorithm
using a test of significance that is developed in this regard.

As motivated above, 
we discard all photometric
information that may be available for a system at hand and use the
kinematic information that is sometimes available, namely velocities
along the line-of-sight (LOS) of individual system members, (such as
stars). We define the galaxy such that the $x_1-x_2$ coordinate system
spans the plane-of-the-sky (POS) and the LOS is along the z-axis.
Thus, our data comprise 1-component $v_3$ values of individual
galactic members and their POS coordinates. We also account for the
errors in $v_3$ measurement. The observed data samples often bear
$\lesssim$ 100 data points. These data are input into the Bayesian
algorithm CHASSIS (Chakrabarty $\&$ Saha 2001, Chakrabarty $\&$
Portegies Zwart 2005).

\section{CHASSIS}
\noindent
CHASSIS helps constrain $\Phi(\cdot)$ and $f(\cdot)$ of galaxies,
given the data $u$ (say). Actually, within CHASSIS, we constrain the
gravitational matter density $\rho({\cdot})$ rather than
$\Phi(\cdot)$, where Poisson equation connects $\rho(\cdot)$ and
$\Phi(\cdot)$ as in:
\begin{equation}
\Phi(\bf{x}) = -4\pi G\nabla^2\rho(\bf{x})
\label{eqn:poisson}
\end{equation}
This helps to avoid problems about negative $\rho(\cdot)$. The
calculation of $\Phi(\cdot)$ from $\rho(\cdot)$ is undertaken at every
step.

Dynamical theory tells us (Dubrovin, Fomenko $\&$ Petrovich, 1990):
\begin{equation}
  f = f[{K_i({\bf{w}})}], \quad {\textrm{where}}\quad \dot{{K_i}}=0\quad \forall i=1,2,3\ldots.
\end{equation}
i.e. ${K_i}$ is an integral of motion. Now, we realise that the size
of the data $u$ is small to moderate and typically bears no more
information other than a single component of velocity. In such a case,
we feel that the data are not sufficient to constrain the extended
forms of $f(\cdot)$ and $\Phi(\cdot)$. In other words, we resort to
making assumptions about $f(\cdot)$ and $\Phi(\cdot)$. 

In fact, we assume $K_i \neq$constant, only for $i=1,2$. Thus,
$f=f[K_1({\bf w}), K_2({\bf w})]$ where $K_1\equiv$energy $E$ of the
galactic particle and $K_2\equiv$ the angular momentum 
or $L$. Here $E=\sum_j v_j^2/2 + \Phi({\bf x})$ and $L=|{\bf
r}\times{\bf v}|$ where ${\bf r}$ is the spherical radius and ${\bf v}$
is the 3-D velocity vector: $\sum_j x_j^2 = r^2$, $v^2=\sum_j
v_j^2$. We also assume radial symmetry in the potential,
i.e. $\Phi({\bf x})= \Phi(r)$. Thus, the particle energy is $E=v^2/2 +
\Phi(r)$. In this geometry, the Poisson equation
(Equation~\ref{eqn:poisson}) is solved numerically by assuming the
mass to be stratified on spherical shells. In other words, the
relevant radial range is discretised and over each radial bin,
$\rho(r)$ is held a constant.

In this background, we seek
\begin{equation}
\Pr(f, \rho|u) \propto \Pr(u| f, \rho)\Pr(f, \rho)
\end{equation}
where the only constraints posed by the priors are physically realistic
requirements of positivity and monotonicity. This owes to the fact that
in general, we do not possess any other prior information that could help
constrain the forms of the sought functions. Thus, 
\begin{eqnarray}
\Pr(f, \rho) &=& \Pr(f)\Pr(\rho), \quad {\textrm{where}}\nonumber \\
\Pr(f) &=& 1 \quad{\textrm {if}}\quad f \geq 0\:\vee {\displaystyle{\frac{\partial f}{\partial E}{\bigg \vert}_{L} \leq 0}}\nonumber \\
{} &=& 0 \quad{\textrm {otherwise}}\nonumber \\
\Pr(\rho) &=& 1 \quad{\textrm {if}}\quad \rho \geq 0\:\vee {\displaystyle{\frac{d\rho}{dr} \leq 0}} \nonumber \\
{} &=& 0 \quad{\textrm {otherwise}}
\end{eqnarray}

We adopt the above priors and estimate the $f(\cdot)$ from which the
data $u$ is most likely to have been drawn, in the estimated
$\Phi(\cdot)$. This is done by iterating towards the most likely set
of $\{f(\cdot), \rho(\cdot)\}$ starting with an arbitrarily chosen
seed. At every step, the current choice of $f(\cdot)$ is projected
into the space of observables (spanned by $x_1, x_2, v_3$); a ready
definition for the likelihood function ${\cal{L}}$ is in terms of the
projection $\eta(\cdot)$ of $f(E, L)$.
\begin{eqnarray}
\eta_i(x_1^i, x_2^i, v_3^i) &=& \displaystyle\int{f[E(v, \Phi(r)), \:L] \:dx_3 dv_1 dv_2} \quad {\textrm {and}} \nonumber \\
{\cal{L}} &=& \displaystyle{\sum_{i=1}^{N_{data}}\log\eta_i}
\label{eqn:integral}
\end{eqnarray}
where $x_1^i, x_2^i, v_3^i$ is the $i^{th}$ data point in the $N_{data}$
sized data sample.

The numerical implementation of a trial $f(\cdot)$ function, at a
trial ($\rho(\cdot)$ or) $\Phi(\cdot)$ is the crucially important
question from the point of view of algorithm design. We do this by
discretising the $E-L$ space and holding $f(\cdot)$ a constant
(=$f_{cell}$) over a given $E-L$ cell. The contribution to
$\eta(\cdot)$, from this $E-L$ cell - defined by say, $E\in[E_1,
E_2], \: L\in[L_1, L_2]$ - is given as:
\begin{eqnarray}
\eta_i^{cell} &=& \displaystyle{f_{cell}\int{dx_3 dv_1 dv_2}} \quad {\textrm {and}} \nonumber \\
\eta_i &=& \displaystyle{\sum_{cell}\eta_i^{cell}}
\end{eqnarray}
Let the integral on the right hand side of the former of these two
equations be $A$. Then, we seek the mapping $A:\longrightarrow
E-L$space.

In order to establish this, first we determine the 2-D area of
intersection between the locii of $E=E_1$, $E=E_2$, $L=L_1$,
$L=L_2$, in the $v_1-v_2$ space. This gives the connection between
the $E-L$ cell at hand and $v_1, v_2$. Mapping $x_3$ to this cell
requires knowledge of the minimum and maximum values of $x_3$ that are
allowed in the cell, given the data point $x_1^i, x_2^i, v_3^i$. This
maximum value is $\sqrt{r_0^2 - (x_1^i)^2 - (x_2^i)^2}$ where $r_0$ is the
solution to: $E_2 = v_3^2/2 + L^2/2r_0^2 + \Phi(r_0)$. The minimum
$x_3$ is 0.

This explains the background to the structures of $\rho(\cdot)$ and
$f(\cdot)$.

\subsection{$\rho$-histogram and $f$-histogram}
\noindent
The representation of $f(E, L)$ over the discretised $E-L$ space,
is akin to a 2-D histogram. Similarly, the $\rho(r)$ structure is
represented as a 1-D histogram. These histograms are updated at the
beginning of every step, while maintaining positivity and
monotonicity. 

The jump distribution we use is discussed below.  If in step $k$, for
$r\in[r_{q-1}, r_q]$ ($\forall q=1,\ldots,N_r$, $r_0$=0), $\rho(r_q) =
\alpha_q^k$, then in the $k+1^{th}$ step:
\begin{equation}
\alpha_{q}^{k+1} = \alpha_{q+1}^k + \displaystyle{(\alpha_{q}^k - \alpha_{q+1}^k)\exp\left(\frac{{\cal R}}{s_1}\right)}
\end{equation}
where ${\cal R}$ is a random number with ${\cal R}\in[-0.5, 0.5]$ and
$s_1$ is an experimentally optimised scale that determines the scale
over which the shape of the $\rho$-histogram is changed. This updating
is done $\forall q$.  Once the shape of the $\rho$-histogram is
updated in this way, the whole structure is scaled by the factor
$\exp({\cal R}/s_2)$ where $s_2$ is another scale.
The $f$-histogram is similarly updated in shape and subsequently normalised.

\subsection{Optimisation}
\noindent
Once the histograms are updated, we project the current $f(\cdot)$
over an $E-L$ cell, into the space of observables, for the $i^{th}$
data point and then sum over all such cells to get $\eta_i$
(Equation~7).  This is done $\forall\:i$, to obtain ${\cal L}$
(Equation~6). The global maxima in ${\cal L}$ is sought by the
Metropolis-Hastings algorithm (Metropolis et. al 1953, Hastings 1970,
Chib $\&$ Greenberg 1995).  Anticipating the likelihood distribution
to be multimodal, we work with highly dispersed seeds to initiate
several chains (Gelman $\&$ Rubin 1992) as well as employ simulated
annealing on a single chain. The latter approach, though perhaps less
obvious, is one that we find very effective in test runs.  

While the
optimisation routine is hard-wired within CHASSIS, the user is allowed
the flexibility to adjust details of the used cooling schedule and
other optimisation parameters such as the scales relevant to the jump
distribution ($s_1$, $s_2$). In this note, it is worth mentioning that 
the current implementation of optimisation is modular, and it is simple
to replace it by a more effcient routine.

\subsection{Required User Input} 
\noindent
The methodology discussed above is incorporated into CHASSIS and all
that user is required to input is the velocity data, the source of
which is independent of CHASSIS. Thus, measured kinematic data,
irrespective of its source, is acceptable, as long as the columns
pertain to the observables $r_p$, $v_3$ and the measured errors in
$v_3$. Here $r_p=\sqrt{x_1^2+x_2^2}$. Besides, the user is allowed to
input details such as the number of bins, bin widths, fraction of data
she wants to perform the run with, the seeds for the sought solutions
and the optimisation related details (see Section~2.2). The user
inputs are advanced via an input file that CHASSIS calls at the
beginning of a run.

\subsection{Assumption of Isotropy}
\noindent
Given the limited data sample, we find that limiting the domain of $f$
to 2-D is not constraining enough in reducing the magnitude of
uncertainties in the estimated solutions to useful levels. Thus, we
resort to imposing the further constraint that $f=f(E)$, i.e. we
assume isotropy to exist in phase space (since $E$ is symmetrical in
$v_j$ and $x_j$, where $j$=1,2,3). Then we resort to (1) justifying or
rejecting our assumption in the data by performing a test of
hypothesis exercise (2)exploring independent measurements that may be
available in the literature to obtain a $\rho$ that is unaffected by
the amount of anisotropy in the data.

\section{Testing for Isotropy - Nonparametric FBST}
\noindent
We test for isotropy in the data, i.e. the null hypothesis $H_0:
{\hat{f}}=\Psi[E(v, r)]$, where ${\hat{f}}$ is the phase space $pdf$
from which the observed data are drawn and $\Psi$ is some function
that manifests phase space isotropy: $\Psi(E)=0\forall E > 0$ and
$\Psi(E) > 0$ otherwise.

This $H_0$ is tested in the data $u$ along the lines of the Fully
Bayesian Significance Test or FBST (Pereira $\&$ Stern 1999; Pereira,
Stern $\&$ Wechsler 2008), except that here, we advance a
nonparametric implementation of the same. Our null is sharp, as is the
requirement for FBST (Madurga, Esteves $\&$ Wechsler 2001). We refer
to our object functions $\{\rho(\cdot), f(\cdot)\}\equiv\theta$ (say);
let $\theta\in\Theta$-space. We assume that $\theta$ is continuous in
the $\Theta$-space. According to FBST, the evidence in favour of $H_0$
is $1- ev$, where:
\begin{eqnarray}
ev &=& 1 - \Pr(\theta \in T\vert u), \quad {\textrm{where}}\quad \nonumber \\
T &=& \{\theta: \Pr(\theta\vert u) > \Pr(\theta^{*}|H_0)\}.
\end{eqnarray}
Here $\theta^{*}$ is the value of $\theta$ which, under the null, 
maximises the posterior $\Pr(\theta\vert u)$.

At the end of every iterative step during a run of CHASSIS, a $\theta$
configuration is identified. The $f(\cdot)$ recovered upon convergence
of the run is indeed a function of $E$ and $E$ only but the true
${\hat{f}}$ is not necessarily so.

Upon convergence of a run of CHASSIS performed with data $u$, we
sample the recovered $\theta$, $N$ times, such that the $i^{th}$
sampling of $\theta$ gives the $i^{th}$ set of observables or $u_i$;
$i=1,2,\ldots,N$. These $N$ data samples are then input into $N$
different new runs of the algorithm. During the run performed with the
data $u_i$, the $j^{th}$ iterative step yields the configuration
$\theta_i^j$ (say), where $\theta_i^j\equiv\{f(\cdot)_i^j,
\rho(\cdot)_i^j\}$. Then $f(\cdot)_i^j$ is isotropic $\forall i,j $, since the
phase space $pdf$ from which $u_i$ is drawn is the recovered $f(E)$,
which by construction, is indeed isotropic.

We scan over all $i, j$ to identify $i^{*},j^{*}$ for which the
posterior is maximised. Thus, $\theta_{i^*}^{j^*}=\theta^{*}$ is
identified. Here $\theta^* \equiv \{f^*(E),\:\rho^*(r)\}$, i.e. functions
recovered at the end of the $j^*$-th step, in a run performed with the
data $u_{i^*}$.

In this nonparametric implementation, $\Pr(\theta \in T\vert u)=X/Y$,
where $X$ is the number of times that a step yields a likelihood in
excess of $\Pr(\theta^{*}|H_0)$ in all the undertaken runs. $Y$ is the
total number of iterative steps in all the runs undertaken.

\section{Application~1 - multistability in galaxies}
\noindent
It is a common practise in astrophysics to employ the measured 1-D
velocity data of suitable galactic members, with the aim of recovering
the total gravitational mass distributions. To examine the viability
of such a practise, we employ the available measured kinematic data of
two distinct classes of galactic members in an example galaxy - as an
aside, these are planetary nebulae (PNe) and globular clusters
(GCs). The data of 164 PNe ($u_P$) are due to Douglas et al (2007) and
that of 30 GCs ($u_G$) are due to Bergond et. al (2006). These sample
sizes are too discrepant to allow for easy interpretation of any
$p$-value based testing of $H_0$ defined in Section~3. Instead, we
resort to the non-parametric FBST discussed above.

\begin{figure}
\centering{
  $\begin{array}{c c}
  \includegraphics[height=.19\textheight]{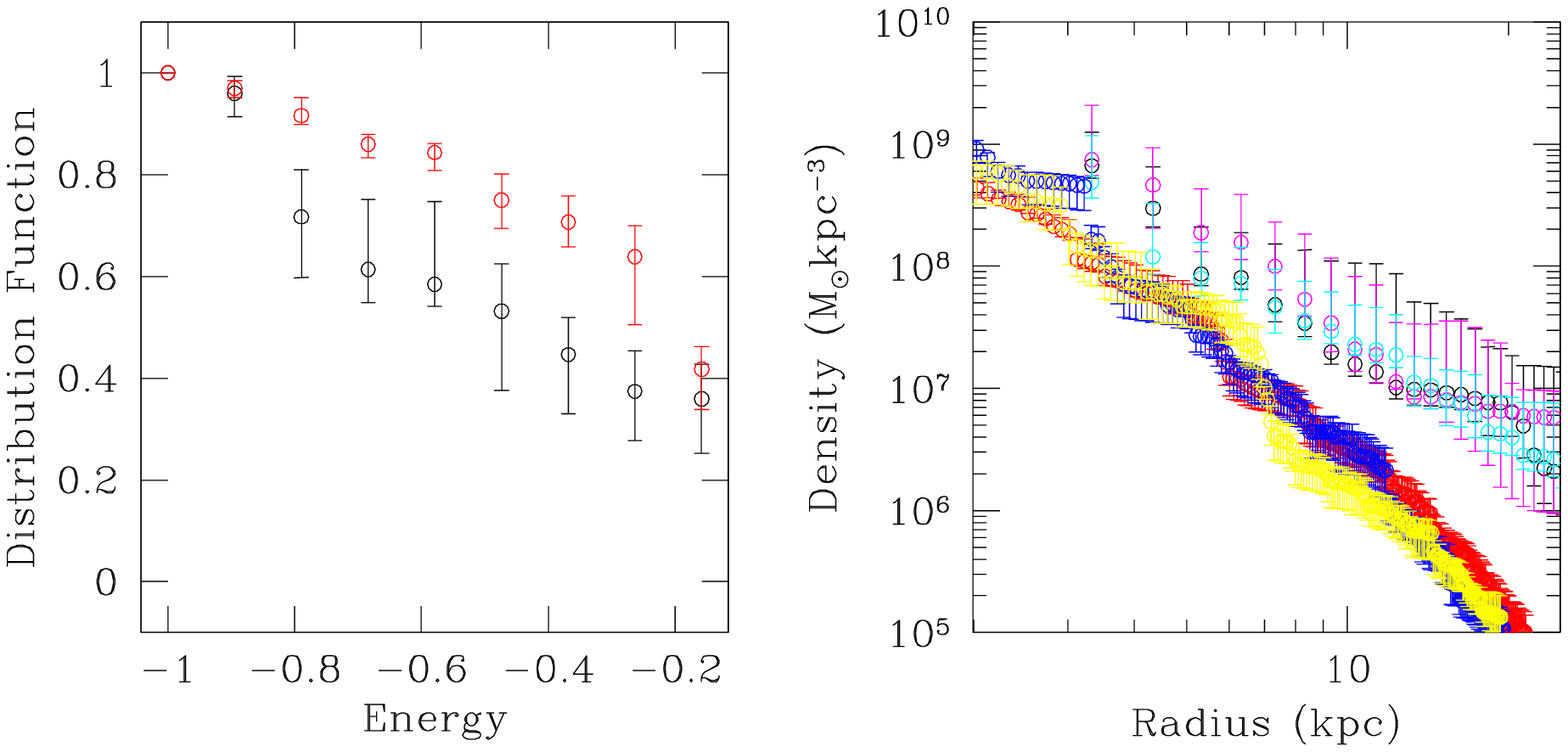} &
  \includegraphics[height=.19\textheight]{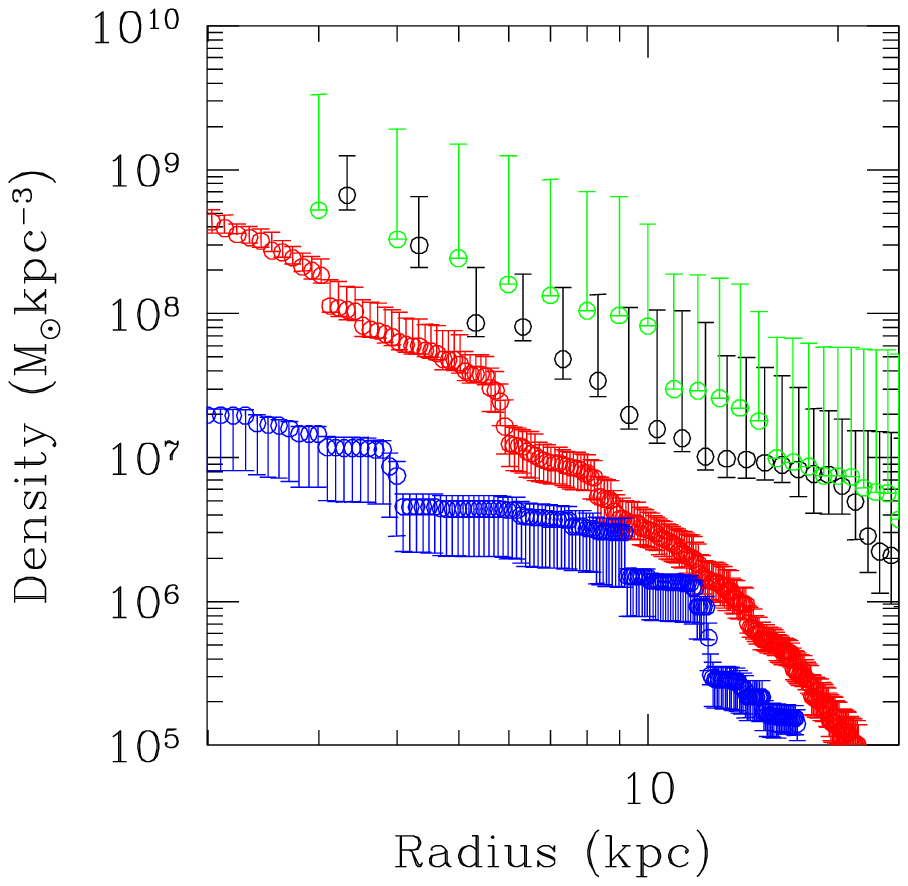} \end{array}$
  \caption{{\it {Left:}} The (normalised) phase space distribution
  functions recovered from two runs of CHASSIS, performed with the GC
  data $u_G$ (in black) and PNe data $u_P$ (in red). The normalisation
  is performed to ensure that $f(E)$=1 for $E$=-1. These $pdf$s are
  recovered under the assumption of isotropy and our implementation of
  nonparametric FBST shows that the profile in black is expected to be
  closer to the true phase space $pdf$ than is the profile in red. The
  errors are $\pm$1-$\sigma$ uncertainties identified on the solution
  by the optimiser. {\it {Middle:}} Gravitational matter density
  distributions from three runs with $u_G$ performed with three
  different seeds (in magenta, cyan and black) and from 3 runs with
  $u_P$ (in yellow, red and blue). The units of density and radius are
  astrophysical. {\it {Right:}} $\rho(r)$ from one of the runs done
  with $u_G$ is shown in black while that with $u_P$ is in red. The
  $\rho^*(r)$ corresponding to the implementation of data $u_G$ is
  shown in green while that corresponding to $u_P$ is in blue. We
  notice that $\rho_G(r)$ is consistent with $\rho^*(r)$ obtained from
  $u_G$ while $\rho^*(r)$ obtained from $u_P$ is significantly lower
  than $\rho_P(r)$.  }}
\label{fig:anchassis}
\end{figure}

The two data sets are input into the {\it{isotropy-assuming}} and
{\it{sphericity-assuming}} CHASSIS.  The $f(E)$ and $\rho(r)$
recovered from the two distinct data are found to be inconsistent with
each other, within error bars (Figure~1). The difference in the
recovered $f(\cdot)$ could arise from different divergences between an
isotropic phase space $pdf$ and the ${\hat{f}}$ from which $u_P$ are
drawn, as compared to that from which $u_G$ are drawn. However, the
distribution of gravitational matter in the galaxy should be uniquely
determined. That such is not our conclusion, prompts us to examine if
the inherent assumption of isotropy is to be blamed. Thus, we test for
our null $H_0$ (defined above in Section~3) in the data $u_P$ and
$u_G$ separately.

The results of our implementation of FBST are shown in Figure~1. We
find that for three different runs done with distinct seeds, given
$u_P$, the average evidence in favour of the null is about 0.60 For
three runs done with different seeds, given $u_G$, the average $1 - ev$ is
0.95. Thus we conclude that the degree of isotropy of the $pdf$ that
$u_G$ is drawn from is higher than that of the $pdf$ that $u_P$ is
drawn from.

We wonder if the difference in the recovered mass distributions be
due to the concluded difference in isotropy in the two data sets? To
understand this, we invoke the peculiarity of CHASSIS that the
algorithm overestimates mass density at all $r$ where phase space
anisotropy prevails (discussed in Section~5). Thus, we expect that
$\rho(r)$ recovered using $u_P$ ($\rho_P(r)$) is more of an
over-estimate compared to the galactic mass density than is the
$\rho(r)$ recovered using $u_G$ ($\rho_G(r)$).

However, as shown in Figure~1, at all $r \gtrsim$ 6 kpc, $\rho_G(r) >
\rho_P(r)$. Therefore, to reconcile the difference between $\rho_P(r)$
and $\rho_G(r)$, the isotropy issue cannot help unless we propose that
$u_G$ is drawn from a more anisotropic $pdf$ than $u_P$. This is of
course not true but its inverse is. Hence we conclude that differences
in $\rho_G(r)$ and $\rho_P(r)$ are intrinsic to the system and not due
to our assumption of isotropy.

Thus we have demonstrated the potential risk in employing kinematic
data of individual galactic members of a particular population type,
to compute the gravitational mass distribution of galaxies. We have
also shown that the galactic phase space is described by at least two
distinct basins of attractions, i.e. the galaxy is multistable, as we
would expect complex systems like galaxies to be.

\section{Application~2 - using the total mass constraint}
\noindent
Motivated by the need to simplify our analysis by reducing the number
of degrees of freedom to the bare minimum, we persist with the
assumption of $f=f(E)$, i.e. phase space is isotropic. We envisage
that when the data have been drawn from an anisotropic $pdf$, the
algorithm will imply erroneous answers. Using physical arguments we
can predict the nature of this error - CHASSIS overestimates $\rho(r)$
at $r$ where phase space anisotropy prevails. 

We implement an independent measure of the total gravitational matter
($M_0$) inside a given radius ($R_E$) in an example galaxy with the
aim of recovering the correct solution for $\rho(r)$ under our
assumption of isotropy, irrespective of the degree of anisotropy of
the $pdf$ from which the input data is drawn. If the system is at a
large distance from us, we cannot get velocity data of individual
galactic members. Then, we can only get projected velocity dispersion
values ($\sigma_p$) at 1 or a few radial locations in the galaxy. In this
background of sparse and incomplete velocity data, we prepare velocity
data samples for inputting to CHASSIS, in the following way.

We select observables ($x_1, x_2, v_3$) from 3 different toy
$f(\cdot)$, two of which are selected to depend on $E$ and $L$ while
the third is isotropic. For the example galaxy, $\sigma_p$ is known at
$r_1$, $r_2$ and $r_3$ say. Then we consider the system to be divided
into 3 anulii with $r\in[0, r_1], (r_1, r_2], (r_2, r_3]$. The error
$\delta$ in the $\sigma_p$ measured at $r_i$ ($i=1, 2, 3$) can be
related to the size of the data sample $N_i$ we intend to draw from
the given annulus, assuming normal error distribution. The phase space
$pdf$ that we choose our samples from are $f_{Gauss},
\:f_{WD}\:f_{Michie}$:
\begin{eqnarray}
f_{Gauss}(E) &=& \displaystyle{\frac{1}{\sqrt{2\pi\sigma^2}}}\displaystyle{exp\left(\frac{-E}{\sigma^2}\right)} \quad E < 0, \nonumber \\
 &=& 0 \quad E > 0, \nonumber \\
f_{WD}(E, L) &=& \displaystyle{\frac{1}{\sqrt{2\pi\sigma^2}}}\displaystyle{\exp\left(-\frac{L^2}{r_a\sigma^2}\right)\exp\left(\frac{-E}{\sigma^2}\right)} \quad E < 0, \nonumber \\
 &=& 0 \quad E > 0. \nonumber \\
f_{Michie}(E, L) &=& \displaystyle{\frac{1}{\sqrt{2\pi\sigma^2}}}\displaystyle{\exp\left(-\frac{L^2}{r_a\sigma^2}\right)\left[\exp\left(\frac{-E}{\sigma^2}\right) - 1 \right]} \quad E < 0, \nonumber \\
 &=& 0 \quad E > 0.
\end{eqnarray}
The samples chosen from $f_{Gauss}, \:f_{WD}\:f_{Michie}$ are
$S_{Gauss}, \:S_{WD}\:S_{Michie}$ (say). Also, to test for the effect
of data from different forms of $f(\cdot)$, we define $E=v^2/2 +
\Phi_{test}(r)$ where we choose $\Phi_{test}(r) \sim
1/\sqrt{r^2 + r_c^2}$. 

\begin{figure}
\centering{
  \includegraphics[height=.2\textheight]{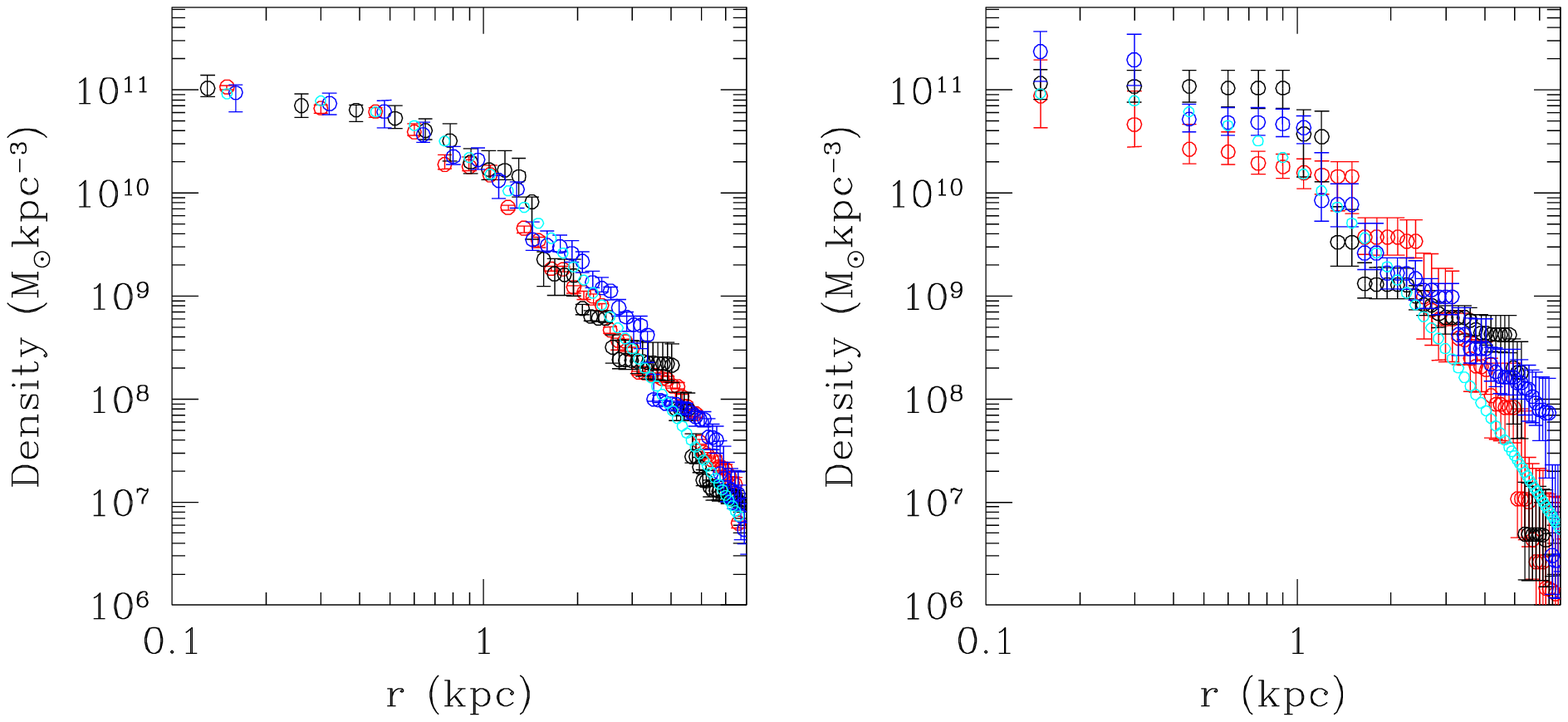} 
  \caption{Gravitational mass density distributions over radius,
obtained from runs done with data drawn from three different phase
space distributions ($f_{Gauss}, \:f_{WD}\:f_{Michie}$) that are
distinguished from each other in terms of the inherent degree of
anisotropy in their forms. $\rho(r)$ estimated from run performed with
data $S_{isotropy}$ is in black, with $S_{WD}$ is in red and with
$S_{Michie}$ is in blue. The $\rho(r)$ implied by the fiduciary test
potential $\Phi_{test}(r)$ is in cyan. The left panel displays results
obtained when the constraint of total mass is included ($M_0 \approx
4.06\times10^{11}$ M$_{\odot}$ within about 8.7 kpc). On the right,
the general inconsistency between the estimated profiles, when the
mass constraint is excluded, is brought out. }}
\label{fig:modsim}
\end{figure}

The total mass constraint is expected to narrow down the 
range of solutions possible; in this sense it acts as a prior
on the solution for $\rho(\cdot)$. 
We test if the chosen data samples, when input into CHASSIS, recover a
$\Phi(r)$ that is concurrent with $\Phi_{test}(r)$ when we
include/exclude the constraint that total mass within $r=R_E$ is $M_0$.
Here $M_0$ is obtained from literature as about 4.06$\times10^{11}$
M$_\odot$, with errors of $\delta M_0 = \pm 0.2\times 10^{11}$
M$_\odot$ and $R_E\approx$8.7 kpc (Koopmans $\&$ Treu 2003). We
incorporate the constraint by adding a penalty function to the
definition of the likelihood; the role of this penalty function is
to penalise any solution 
that implies a total
gravitational mass within $R_E$ ($M_c(R_E)$) different from
$M_0$. This penalty function is $\alpha\vert M_c(R_E) -
M_0\vert/2\delta M_0$. Here $\alpha$ is a flag, designed to include
or exclude the constraint from the definition of the likelihood, depending
on whether it is 1 or 0 respectively.

The results of conducting our test runs are shown in Figure~2. We
find that when the constraint is included, the $\rho(r)$ profiles are
consistent with each other within errors, irrespective of the
anisotropy in the data used to obtain this profile.
In absence of guidance from the constraint, anisotropy affects results.

\section{Summary}
\noindent
Here we have discussed the novel nonparametric Bayesian algorithm
CHASSIS that estimates the most likely phase space $pdf$ ($f(\cdot)$)
from which an observed sample of 1-D component of velocities of
individual galactic members is drawn, at the most likely gravitational
potential $\Phi(\cdot)$ of the galaxy.  The main purpose of this paper
is to bring out the fact that CHASSIS is a {\it robust and viable
  algorithm} that can be implemented to extract the all-important
gravitational matter density distribution in distant galaxies, even
within the domain of very sparse and incomplete data. Given the dearth
of measurements in these systems, it is prudent to work with a small
number of degrees of freedom. With this in mind, one version of
CHASSIS has been designed to work under the purview of the assumption
of isotropy in phase space, by which we imply an $f(\cdot)$ that is
symmetric in the 3 velocity and 3 spatial coordinates. In another
version, the assumption of isotropy is not made and CHASSIS is made to
work with a greater number of dof (Chakrabarty $\&$ Saha, under
preparation).

We offer independent means of tackling the obvious fallout of the
assumption of isotropy, when invoked; phase space $pdf$s from which
realistic data are drawn, will not be isotropic in general, leading to
spurious solutions for $f(\cdot)$ and $\Phi(\cdot)$.  Firstly, a
robust test of hypothesis is developed that tests the assumption of
isotropy, given the data. This test is modelled after Pereira $\&$
Stern's (1999) FBST and is designed to work in the nonparametric
context. Data available in astrophysical literature are implemented to
conclude that kinematic data of distinct galactic member populations
will in general provide distinct gravitational mass distributions. The
multistability of the example galaxy is demonstrated and the folly of
this mode of mass determination is indicated. Secondly, we demonstrate
that the usage of information about total gravitational mass,
available in the literature, can constrain the sought gravitational
matter density distribution, irrespective of anisotropy in the
data. Such a constraint basically supplements for the uninformative
priors that we use within CHASSIS.
 
Currently, we are working on the establishment of a critical value of
the evidence value in favour of a given null. This will enable the
quantified judgement of when to reject or accept the null, given the
data. At the moment, our implementation of FBST only allows for a
comparative judgement.

\end{document}